\documentclass[reqno]{amsart}
\usepackage{hyperref}

\usepackage{amsmath}
  \usepackage{paralist}
  \usepackage{graphicx}
  \usepackage{amssymb}

\numberwithin{equation}{section}

\def \beq{\begin{equation}}
\def \eeq{\end{equation}}

\usepackage{color}

\usepackage{setspace}

\begin{document}
\title[elastodynamics model]{Wave Interaction For A System In Elastodynamics With Initial Data Lying On The Level Set Of One Of The Riemann Invariants}
\author[Kayyunnapara Divya Joseph]{ Kayyunnapara Divya Joseph$^1$ \\
 PhD, Visvesvaraya Technological University, \\
  Belagavi,  Karnataka 590018, India 
}
\email{divediv3@gmail.com}

\maketitle

\begin{abstract}
This paper is concerned with the study of interaction  of waves originating from the Riemann problem centred at two different points for  a  system of equations   modelling propagation of elastic
waves. The system consists of two equations for  $(u,\sigma)$, where, u is the velocity and $\sigma$ is the stress and is strictly hyperbolic and nonconservative.  Study of interaction of waves is one of the  most  important
steps in the construction of global solution with initial data in the space of functions of bounded variation  using approximation procedure like the Glimm's scheme. This amounts to constructing a solution  with initial data consisting of three states $(u_L, \sigma_L), (u_m, \sigma_m),$ and  $(u_R, \sigma_R)$. Usually this analysis is done for the states which are in a small neighbourhood of a fixed state.  Here we get explicit formula for the solution of the system  
when the  data lies in the level sets of Riemann invariants.  The speciality of the work is that   we donot assume smallness conditions on the initial data.  


\end{abstract}

\section{Introduction}

In this paper we study a $2 \times 2$  system of equations modelling the propagation of elastic
waves, namely,
\begin{equation}
\begin{gathered}
\frac{\partial u}{\partial t }+ u \frac{\partial u}{\partial x} - \frac{\partial{\sigma}}{\partial x} =
0,\\
\frac{\partial{\sigma}}{\partial t} + u \frac{\partial{\sigma}}{\partial x} - k^2
\frac{\partial{u}}{\partial x} = 0.
\end{gathered}
\label{e1.1}
\end{equation}
 This
system was derived by Colombeau and Leroux \cite{c1} 
from the one dimensional model of the
equation of elasticity when density is
assumed to remain close to a constant.
Here $u$ is the velocity, $\sigma$ is the stress and $k>0$
is the speed
of propagation of the elastic waves. 

 It is well known that smooth global in time solutions do not exist even if the initial data
is smooth. This system is nonconservative, the product $u \sigma_x$ appearing in the equation does not make sense in the sense of distributions and classical Lax- Glimm theory of conservation laws does not apply. Cauret et al.\cite{c2} and Joseph \cite{j1} constructed solution of the  initial value problem in the class of Colombeau algebra. There is an extensive amount of lierature where  the initial value problem to this system \eqref{e1.1} has been studied in the class of generalized functions. To our knowledge, this work is the first of its kind to study the structure and interaction of waves when the initial data consists of three constant states lying in level sets of one of the Riemann invariants. The study of the case when the initial data lies in the other regions, is more complex and is being studied on. The significance of the work is, in  not having any smallness assumptions on the data. Study of Riemann problems and Interaction of waves is important in  the construction of solutions for the initial value problem in the class of functions of bounded variation using  Glimm scheme.  Rigorous study of initial value problem in the space of functions of bounded variation, for nonconservative system was initiated by Volpert \cite{v1} and developed further  by LeFloch \cite{le1} , LeFloch, Dal Maso and Murat \cite{c3}and LeFloch and Liu \cite{le3}.  Any
discussion of well-posedness of solution for nonconservative system should be based on a given nonconservative product in addition to admissibility criterion for shock discontinuities. 

Existence of solutions for initial value problem for  general non conservative systems, with  initial data of small total variation was proved by LeFloch \cite{le1} and LeFloch and Liu \cite{le3}, using  Glimm scheme \cite{x1} . Here one of the main steps is to solve the
Riemann problem, i.e. with initial data consisting of two constant states separated by a discontinuity. Glimm scheme involves piecing together of solution of Riemann problems with discontinuity at different points. Initially at $t=0$, the initial data is approximated as constants on mesh length $2\Delta x$ in the space variable and solve a family of Riemann problems with discontinuities at $2 m \Delta x, m=0, \pm 1, \pm 2,... $. The  mesh height  $\Delta t$ is chosen such that,  nearby Riemann problems do not interact and  constant states for the  the next step, for the Riemann problems are found  by a Random choice method  to get approximations in the next time step as solution of Riemann problems. Here it is important to understand  
the interaction of nearby Riemann problems.  In other words initial data having two discontinuities and taking three constant states of the following form : 
\begin{equation}
( u(x,0),\sigma (x, 0) )= \begin{cases}
(u_L,\sigma_L),&\text{if } x< x_0,\\
 (u_m,\sigma_m),&\text{if } x_0 < x < x_1\\
(u_R,\sigma_R) ,&\text{if } x > x_1
\end{cases}
\label{e1.2}
\end{equation}
For small data, this study is done in earlier mentioned papers. Here we write down explicitly the solution of this interaction problem for large data lying on the level surface of the Riemann invariants. We use the Volpert product \cite{v1} without any smallness condition on the data. First we recall
some known facts about the Riemann problem for
\eqref{e1.1}.
Riemann problem with no smallness condition on the data was solved in \cite{j2}. Initial Boundary value problem for \eqref{e1.1} is studied in \cite{d1, d2}, for special types of initial data. In this paper we study the interaction
problem \eqref{e1.1} and \eqref{e1.2} with no smallness condition on the states $(u_L,\sigma_L), (u_m,\sigma_m)$ and $(u_R,\sigma_R)$
We write \eqref{e1.1} as the system
\begin{equation*}
\begin{bmatrix}
u\\  
\sigma
\end{bmatrix} 
_t +
\begin{bmatrix}
u & -1\\
-k^2 & u
\end{bmatrix}
\begin{bmatrix}
u \\
\sigma
\end{bmatrix} 
_x
=0
\end{equation*}
The eigevalues of the matrix
\begin{equation*}
A(u,\sigma) =
\begin{bmatrix}
u & -1\\
-k^2 & u
\end{bmatrix}
\end{equation*}
 are called the
characteristic speeds of system \eqref{e1.1}. A simple
computation
shows that, the equation for the eigenvalues are given
by
\begin{equation*}
\lambda_1 = u-k ,\,\,\,\lambda_2 = u + k.
\end{equation*} 
Call $\ r_1 $ and $\ r_2 $ right eigenvectors
corresponding to $ \lambda_1 $ and $ \lambda_2 $
respectively. An
 easy computation shows
 \begin{equation*}
r_1 = r_1(u,\sigma)=(1,k),\,\,\,r_2= r_2 (u,
\sigma)=(1,-k)
\end{equation*}
We say $w_i :R^2 \rightarrow R$ is called $i$-Riemann
invariant if $Dw_i(u,v).r_i(u,\sigma)=0$.
This means that $w(u,\sigma)$ is a $1-$ Riemann invariant
if
$w_u + k w_{\sigma} = 0$ and $2-$ Riemann invariant if
$w_u -k w_{\sigma} = 0$.
Now  $ w_u +k w_{\sigma} = 0 $ gives
$w (u,\sigma)=\phi_1(\sigma - k u) $, where $\phi_1 :R^1
\rightarrow R^1$ is arbitrary smooth function,
whereas $w_{\sigma} - k w_u = 0$ gives
$w(u,\sigma)=\phi_2(\sigma + k u)$, $\phi_2:R^1\rightarrow
R^1$.
This calculation showed that, the functions $w_1$ and $w_2$ given by 
\begin{equation*}
w_1(u,\sigma) =\sigma - k u,\,\,\,
w_2(u,\sigma) =\sigma + k u 
\end{equation*}
are $1-$ Riemann invariant and $2-$Riemann invariant
respectively. Further it was proved in \cite{d1} that $w_1$ is constant along $2-$ characteristic curve
and $w_2$ is constant along $1-$ characteristic curve.
\\

\section{Weak formulation and Rankine-Hugoniot condition}
To explain the weak formulation of solution we introduce the
nonconservative product of Volpert \cite{v1} in sense of measures.
For $w$  a function of bounded variation in $\Omega$ open subset of $R\times (0,\infty)$, we write
\[
\Omega=S_c \cup S_j \cup S_n
\]
where $S_c,$ $S_j$ are the points of approximate continuity of $w$ and  the
points of approximate jump of $w$ respectively and $S_n$ is a set of one dimensional
Hausdorff-measure zero. For a.e.  point $(x,t) \in S_j$, with respect to the one dimensional Hausdrorff measue there exists a unit vector $\nu=(\nu_x,\nu_t)$  unit normal to $S_j$ at $(x,t)$ and $w(x-0,t)$ and
$w(x+0,t)$ called the left and right values of $w(x,t)$.   For any
continuous function $g :R^1 \to R^1$, and a function $v$ of bouned variation in $\Omega$ the Volpert product
$g(w)v_x$ is defined as a Borel measure in the following manner.
Consider the averaged superposition of $g(w)$  (see Volpert \cite{v1})
\begin{equation}
\overline{g(w)}(x,t) = \begin{cases}
 g(w(x,t)),&\text{if } (x,t) \in S_c,\\
  \int_0^1 g((1-\alpha)(u(x-,t)+\alpha u(x+,t))d\alpha,
  &\text{if }(x,t) \in S_j
\end{cases}
\label{m2.1}
\end{equation}
and the associated measures $[g(w)v_x]$ and 
$[g(w)v_t]$ defined by
\begin{equation}
[g(w)v_x](A)=\int_{A}\overline{g(w)}(x,t)v_x,\,\, [g(w)v_t](A)=\int_{A}\overline{g(w)}(x,t)v_t,
\label{m2.2}
\end{equation}
where $A$ is a Borel measurable subset of $S_c$ and
\begin{equation}
\begin{aligned}
&[g(w)v_x](\{(x,t)\})=\overline{g(w)}(x,t)(v(x+0,t)-v(x-0,t))\nu_x,\\
&[g(w)v_t](\{(x,t)\})=\overline{g(w)}(x,t)(v(x+0,t)-v(x-0,t))\nu_t
\end{aligned}
\label{m2.3}
\end{equation}
provided $(x,t) \in S_j$.

The system of  equations in \eqref{e1.1} is understood as

\begin{equation}
\begin{gathered}
\frac{\partial u}{\partial t }+\bar{ u} \frac{\partial u}{\partial x} - \frac{\partial{\sigma}}{\partial x} =
0,\\
\frac{\partial{\sigma}}{\partial t} + \bar{u} \frac{\partial{\sigma}}{\partial x} - k^2
\frac{\partial{u}}{\partial x} = 0.
\end{gathered}
\label{m2.4}
\end{equation}
in the sense of measures, where $\bar{ u}$ is defined by \eqref{m2.1}, \eqref{m2.2} with $g(u) =u$. The Rankine Hugoniot condition is the statement that this measure is zero on the jump set.

 To get the Rankine Hugoniot condition, we compute this measure on points of approximate jump $S_j$ and set it equal to zero.
Consider a surface of discontinuity of $u$ given by $\gamma(x,t)=x-\phi(t)=0$. Then the normal at a point $(x,t)$ on the surface is given by 
\begin{equation}
(\nu_t,\nu_x)=\frac{(\gamma_t,\gamma_x)}{|(\gamma_t,\gamma_x)|}= \frac{1}{\sqrt{1+\phi'(t)^2}}(-\phi'(t),1).
\label{m2.5}
\end{equation}
 From the definition of Volpert product, along the curve of discontinuity $x=\phi(t)$, we have using \eqref{m2.3} and \eqref{m2.4}, 
\begin{equation}
\begin{aligned}
&\nu_t (u(\phi(t)+)-u(\phi(t)-)+\frac{u(\phi(t))+u(\phi(t)-)}{2}(u(\phi(t)+-u(\phi(t)-)\nu_x\\
&-(\sigma(\phi(t)+)-\sigma(\phi(t)-))\nu_x=0\\
&\nu_t (\sigma(\phi(t)+)-\sigma(\phi(t)-)+\frac{u(\phi(t))+u(\phi(t)-)}{2}(\sigma(\phi(t)+)-\sigma(\phi(t)-)\nu_x\\
&-k^2(u(\phi(t)+)-u(\phi(t)-))\nu_x=0
\end{aligned}
\label{m2.6}
\end{equation}

Using \eqref{m2.5} in \eqref{m2.6}, at the point of discontinuity $(\phi(t),t)$, with $u_\pm=u(\phi(t)\pm)$ and $\sigma_\pm=\sigma(\phi(t)\pm)$, $s=\phi'$ we get the Rankine Hugoniot condition \eqref{m2.7},
\begin{equation}
 \begin{gathered}
-s(u_{+}-u_{-})+\frac{u_{+}^2 - u_{-}^2}{2}
-(\sigma_{+}-\sigma_{-})=0\\
-s(\sigma_{+}-\sigma_{-})+\frac{u_{+}+u_{-}}{2}(\sigma_{+}-\sigma_{-})-k^2(u_{+}-u_{-})
\end{gathered}
\label{m2.7}
\end{equation}

\section{The Riemann problem }
 First we recall
some known facts about the Riemann problem for
\eqref{e1.1}, see \cite{j2} for details. Here the initial data takes the form
\begin{equation*}
(u(x,0),\sigma(x,0))  = \begin{cases}
(u_L,\sigma_L),&\text{if } x< 0,\\
(u_R,\sigma_R) ,&\text{if } x > 0
\end{cases}
\end{equation*}
and there is no smallness conditions on the data. A 1-  rarefaction wave connecting $( u_R, \sigma_R )$ on the right and $(u_L, \sigma_L)$ on the left,  is a weak solution of \eqref{e1.1} of the form
\begin{equation*}
(u(x,t),\sigma(x,t)) = \begin{cases}
(u_L,\sigma_L),&\text{if } x < (u_L - k) t,\\
( \frac{x}{t} +  k, k \frac{x}{t} + \sigma_L - k ( u_L - k) )  &\text{if }  (u_L - k) t < x < (u_R - k) t,\\
 (u_R,\sigma_R),&\text{if } x >  (u_R - k) t. 
\end{cases}
\end{equation*}
A 2-  rarefaction wave connecting $( u_R, \sigma_R )$ on the right and $(u_L, \sigma_L)$ on the left,  is a weak solution of \eqref{e1.1}  of the form
\begin{equation*}
(u(x,t),\sigma(x,t)) = \begin{cases}
(u_L,\sigma_L),&\text{if } x < (u_L - k) t,\\
( \frac{x}{t} -  k, - k \frac{x}{t} + \sigma_L + k ( u_L + k) )                      
    &\text{if }  (u_L - k) t < x < (u_R - k) t,\\
 (u_R,\sigma_R),&\text{if } x >  (u_R - k) t. 
\end{cases}
\end{equation*}
A j -  shock wave  connecting $( u_R, \sigma_R )$ on the right and $(u_L, \sigma_L)$ on the left,  is a weak solution of \eqref{e1.1} with
speed $s$  of the form
\begin{equation*}
(u(x,t),\sigma(x,t)) = \begin{cases}
(u_L,\sigma_L),&\text{if } x < {s_j} t,\\
 (u_R,\sigma_R),&\text{if } x > {s_j} t. 
\end{cases}
\end{equation*}

In \cite{j2}, the Riemann problem was solved using Volpert
product and with Lax's admissibility for shocks. It was
shown that
corresponding to each characteristic family $\lambda_j,
j=1,2$ we can define shock waves and rarefaction waves.
Fix a state $(u_{-},\sigma_{-})$
the set of states $(u_{+},\sigma_{+})$ which can be
connected by a single $j$- shock waves is a straight line
called
j-shock curve and is denoted by $S_j$ and the states which
can be connected by a single $j$-rarefaction wave is a
straight line is called $j$ rarefaction curve
and is denoted by $R_j$ and
\begin{equation*}
\begin{gathered}
R_1(u_{-},\sigma_{-}): \sigma=\sigma_{-}+k(u-u_{-}),
u>u_{-}\\
S_1(u_{-},\sigma_{-}): \sigma=\sigma_{-}+k(u-u_{-}),
u<u_{-}\\
R_2(u_{-},\sigma_{-}): \sigma=\sigma_{-}-k(u-u_{-}),
u>u_{-}\\
S_2(u_{-},\sigma_{-}): \sigma=\sigma_{-}-k(u-u_{-}),
u<u_{-}\\
\end{gathered}
\end{equation*}
Further the $j$- shock speed $s_j$is given by 
\begin{equation*}
 s_j =\frac{u_{+}+u_{-}}{2} +(-1)^j k, \,\,\, j=1,2.
\end{equation*}
The Lax entropy condition requires that the $j$- shock
satisfies inequality
\begin{equation*}
\lambda_j(u_{+},\sigma_{+}) \leq s_j \leq
\lambda_j(u_{-},\sigma_{-}).
\end{equation*}
This curves fill in the $u-\sigma$ plane and the Riemann
problem can be solved uniquely for arbitrary initial
states
$(u_{-},\sigma_{-})$ and $(u_{_+},\sigma_{+})$
in the class of self similar functions consisting of
solutions of shock waves and rarefaction waves
separated by constant states. These constant states are
obtained from the shock curves and rarefaction curves
corresponding to the two families of the characteristic
fields.

\section{Interaction of waves} In this section we consider interaction of waves originating at two different points. More specifically we consider solutions for the initial value problem to \eqref{e1.1},
with initial data of the form \eqref{e1.2}. We need to consider different cases depending the position of  $(u_m,\sigma_m)$  relative to $(u_L,\sigma_L)$ and the position of $(u_R,\sigma_R)$ relative to $(u_L,\sigma_L)$ and $(u_m,\sigma_m)$. Fix  $(u_L,\sigma_L)$, and consider the wave  curves 
\begin{equation*}
\begin{gathered}
R_1(u_{L},\sigma_{L})= \{(u,\sigma) \sigma=\sigma_{L}+k(u-u_{L}),
u>u_{L}\}\\
S_1(u_{L},\sigma_{L})= \{(u,\sigma): \sigma=\sigma_{L}+k(u-u_{L}),
u<u_{L}\}\\
R_2(u_{L},\sigma_{L})=\{(u,\sigma): \sigma=\sigma_{L}-k(u-u_{L}),
u>u_{L}\}\\
S_2(u_{L},\sigma_{L})=\{(u,\sigma): \sigma=\sigma_{L}-k(u-u_{L}).
u<u_{L}\}\\
\end{gathered}
\end{equation*}
passing through $(u_L,\sigma_L)$ . These curves divide the $u-\sigma$ plane into 4 regions. 
 Let  $\Gamma_1$ denote the region between $R_1(u_L,\sigma_L)$ and  $R_2(u_L,\sigma_L)$,  $\Gamma_2$ be the region between $R_2(u_L,\sigma_L)$ and  $S_1(u_L,\sigma_L)$, $\Gamma_3$ be the region between $S_1(u_L,\sigma_L)$ and  $S_2(u_L,\sigma_L)$, $\Gamma_4$ be the region between $S_2(u_L,\sigma_L)$ and  $R_1(u_L,\sigma_L)$.
 We now consider the case when $( u_L, \sigma_L ), ( u_m, \sigma_m ), ( u_R, \sigma_R ),$ lie in each of the eight regions,  i.e when $( u_m, \sigma_m ), ( u_R, \sigma_R )$ lies in each of the regions $R_i(u_L,\sigma_L)$ and  $S_i(u_L,\sigma_L), i= 1, 2$ and $\Gamma_i, i=1, 2, 3, 4.$  In order to do this, we first write explicitly the situation when  $( u_L, \sigma_L ), ( u_m, \sigma_m ), ( u_R, \sigma_R ),$ all lie on the level set of the same Riemann Invariants, i.e either on the 1-  Riemann Invariant or 2-  Riemann Invariant. The cases when they do not all lie on the same  Riemann Invariants, has to be dealt with separately. 
\\

{ \bf Theorem : } If $(u_m,\sigma_m), (u_L,\sigma_L), (u_R,\sigma_R)$ lie on the level set of the same Riemann Invariant then there exist a global in time weak solution which can be explicitly written.  \\

{ \bf Proof : } First we consider $u_m, u_L, u_R$ to lie on the level set of the second Riemann Invariant. \\
Case 1: When $(u_m, \sigma_m) \in R_2( (u_L, \sigma_L)  )$ \\
Clearly $u_L$ is connected  to $u_m$ by a $2-$ rarefaction as seen in Figure 1 below. \\
Subcase 1:   If $(u_R, \sigma_R) \in S_2( (u_L, \sigma_L)  ) \bigcup \{ R_2(u_L, \sigma_L), u_L < u_R< u_m \} :$ \\
 then, $u_m$ is connected  to $u_R$ by a $2-$ shock.  Now, since $u_m> u_R$ we have $u_m + k > \frac{u_R +u_m}{2} + k.$ So  the lines  $x=  (u_m + k)t +x_0$and $ x=  (\frac{u_R +u_m}{2} + k) t +x_1 $ intersect at $( x_2, t_2).$ The solution is continued with shock curve  $x=\beta_1(t), \beta_1(t_2) =x_2$, for $t > t_2,$ is given by, 
\begin{equation}
\begin{aligned}
\frac{dx}{dt} &= \frac{u_R - k + \frac{x - x_0}{t} }{2}- \frac{\sigma_R - ( - k \frac{x - x_0}{t} + \sigma_L + k (u_L + k)  )}{u_R + k - \frac{x - x_0}{t} }, \\ 
\frac{dx}{dt} &= \frac{u_R - k + \frac{x - x_0}{t} }{2}- k^2 ( \frac{u_R + k - \frac{x - x_0}{t} }{ \sigma_R - ( - k \frac{x - x_0}{t} + \sigma_L + k (u_L + k)  ) } ).
\label{4.1}
\end{aligned}
 \end{equation}
\begin{figure}[!hbtp]
\includegraphics[width=12cm,height=4cm]{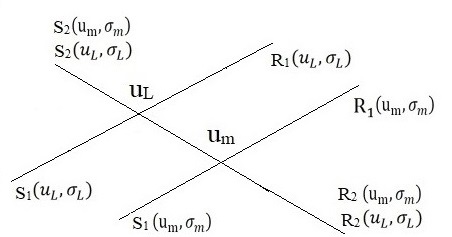}
\caption{ When $(u_m, \sigma_m) \in R_2( (u_L, \sigma_L)  )$}
\end{figure} 
Since  $(u_R, \sigma_R) \in S_2(u_L, \sigma_L)$,  $\sigma_R - \sigma_L - k (u_L - u_R)=0,$. Using this in \eqref{4.1},  the two equations become a single equation and hence the shock curve is given by
\begin{equation*}
\frac{dx}{dt} = \frac{u_R + k }{2} +  \frac{x - x_0}{2 t}, \,\,\, x(t_2) = x_2.
 \end{equation*}
\begin{figure}[!hbtp]
\includegraphics[width=14cm,height=4cm]{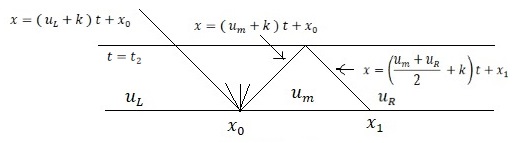}
\caption{ For Case 1,  Subcase 1}
\end{figure} 
Solving this we get
$x= \beta_1(t) = (u_R + k) t + c t^{\frac{1}{2}} +x_0$ where $c= (x_2 - x_0) {t_2}^{- \frac{1}{2}} - (u_R + k) {t_2}^{\frac{1}{2}}.$ Here there are two cases according to whether the curves $x= \beta_1(t)$ and $x= (u_L  + k) t +x_0,$  intersect or not. If the the two curves do not intersect the solution becomes,
\begin{equation}
\begin{gathered}
(u(x,t), \sigma(x, t))  \\
= \begin{cases}
(u_L, \sigma_L),  & \text{if }  x <(u_L + k)t +x_0,  \forall t ,\\
(\frac{x - x_0}{t} - k, - k \frac{x - x_0}{t} + \sigma_L + k (u_L + k) ) , & \text{if }  (u_L + k)t +x_0 < x  < (u_m + k)t  +x_0, \\ & \,\,\, \,\,\, \text{and } t \leq t_2,\\
(u_m, \sigma_m),   & \text{if  }  (u_m + k)t +x_0  <x < (\frac{u_R +u_m}{2}   + k) t +x_1, 
\\ & \,\,\, \,\,\, \text{and }  t \leq t_2, \\
 (u_R, \sigma_R),    & \text{if  } x> (\frac{u_R +u_m}{2} + k) t +x_1,  t \leq t_2, \\
( \frac{x - x_0}{t} - k, - k \frac{x - x_0}{t} + \sigma_L + k (u_L + k)),   & \text{if  } (u_L + k)t +x_0<x<\beta_1(t),  t> t_2,\\
  (u_R, \sigma_R) ,     & \text{if  }   x>\beta_1(t),  t> t_2.
\end{cases}
\label{4.2}
\end{gathered}
 \end{equation}
 If the curves $x= \beta_1(t)$ and $x= (u_L  + k) t +x_0,$  intersect at $(x_3,t_3)$ and if $(u_R, \sigma_R) \in S_2(u_L, \sigma_L),$  then for $t>t_3$, solution is a $2-$ shock with left value $(u_L,\sigma_L)$ and right value $(u_R,\sigma_R)$. So,  if the curves intersect and $(u_R, \sigma_R) \in S_2(u_L, \sigma_L)$  then the solution is, 
\begin{equation}
\begin{gathered}
(u(x,t), \sigma(x, t)) \\
= \begin{cases}
(u_L, \sigma_L),  & \text{if }  x <(u_L + k)t +x_0,  0 \leq t \leq t_2 ,\\
(\frac{x - x_0}{t} - k, - k \frac{x - x_0}{t} + \sigma_L + k (u_L + k) ) , & \text{if }  (u_L + k)t +x_0 < x  < (u_m + k)t  +x_0, 
\\ & \,\,\, \,\,\, \text{and } t \leq t_2,\\
(u_m, \sigma_m),   & \text{if  }  (u_m + k)t +x_0  <x < (\frac{u_R +u_m}{2}  + k) t +x_1, 
\\ & \,\,\, \,\,\, \text{and }  t \leq t_2, \\
 (u_R, \sigma_R),    & \text{if  } x> (\frac{u_R +u_m}{2} + k) t +x_1,  t \leq t_2, \\
( \frac{x - x_0}{t} - k, - k \frac{x - x_0}{t} + \sigma_L + k (u_L + k)),   & \text{if  } (u_L + k)t +x_0<x<\beta_1(t),  t_3> t> t_2,\\
 (u_R, \sigma_R) ,     & \text{if  }   x>\beta_1(t), t_3> t> t_2\\
(u_L,\sigma_L) , & \text{if  } x<(\frac{u_L +u_R}{2} +k)(t-t_3) +x_3, t> t_3,\\
(u_R,\sigma_R) , & \text{if  } x>(\frac{u_L +u_R}{2} +k)(t-t_3, )+x_3,   t> t_3.
\end{cases}
\label{4.3}
\end{gathered}
 \end{equation}
 If the curves $x= \beta_1(t)$ and $x= (u_L  + k) t +x_0,$  intersect at $(x_3,t_3)$ and if $(u_R, \sigma_R) \in R_2(u_L, \sigma_L), u_L < u_R< u_m,$  then for $t>t_3$, solution is a $2-$ rarefaction with left value $(u_L,\sigma_L)$ and right value $(u_R,\sigma_R)$. So,  if the curves intersect and 
$(u_R, \sigma_R) \in R_2(u_L, \sigma_L), u_L < u_R< u_m,$  then the solution is,
\begin{equation}
\begin{gathered}
(u(x,t), \sigma(x, t))\\
= \begin{cases}
(u_L, \sigma_L),  & \text{if }  x <(u_L + k)t +x_0,  0 \leq t \leq t_2,\\
(\frac{x - x_0}{t} - k, - k \frac{x - x_0}{t} + \sigma_L + k (u_L + k) ) , & \text{if }  (u_L + k)t +x_0 < x  < (u_m + k)t  +x_0,
\\ & \,\,\, \,\,\, \text{and } t \leq t_2,\\
(u_m, \sigma_m),   & \text{if  }  (u_m + k)t +x_0  <x < (\frac{u_R +u_m}{2}   + k) t +x_1, 
\\ & \,\,\, \,\,\, \text{and }  t \leq t_2, \\
 (u_R, \sigma_R),    & \text{if  } x> (\frac{u_R +u_m}{2} + k) t +x_1,  t \leq t_2, \\
( \frac{x - x_0}{t} - k, - k \frac{x - x_0}{t} + \sigma_L + k (u_L + k)),   & \text{if  } (u_L + k)t +x_0<x<\beta_1(t),  t_3> t> t_2,\\
 (u_R, \sigma_R) ,     & \text{if  }   x>\beta_1(t), t_3> t> t_2\\
(u_L,\sigma_L) , & \text{if  } x<(u_L + k)(t- t_3) +x_3, t> t_3,\\
(\frac{x - x_0}{t} - k, - k \frac{x - x_0}{t} + \sigma_L + k (u_L + k) ), & \text{if  } (u_L + k)(t- t_3) +x_3< x<  (u_R + k)(t- t_3) 
\\ & \,\,\, \,\,\,  +x_3, t> t_3, \\
 (u_R,\sigma_R)  & \text{if  }  x>  (u_R + k)(t- t_3) +x_3, t> t_3,
\end{cases}
\label{4.4}
\end{gathered}
 \end{equation}
see figure 2. \\

Subcase 2:   If $(u_R, \sigma_R) \in R_2( (u_L, \sigma_L)  ), \,\,\, u_L< u_m < u_R:$ \\
\begin{figure}[!hbtp]
\includegraphics[width=14cm,height=4cm]{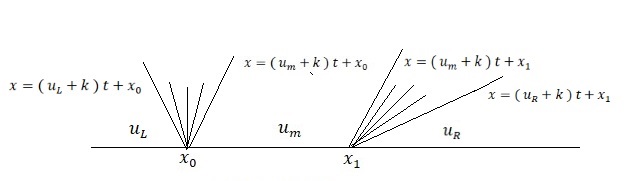}
\caption{ For Case 1,  Subcase 2}
\end{figure} 
Since $(u_L + k)t <  (u_m + k)t <   (u_R + k)t$ and $x_0<x_1,$ there is no interaction and  $(u_L,  \sigma_L)$ is connected to $(u_m,  \sigma_m)$ and $(u_m,  \sigma_m)$ to $(u_R,  \sigma_R)$ by a $2$- rarefaction. The solution here becomes,
\begin{equation}
\begin{gathered}
( u(x,t), \sigma(x, t) ) \\
= \begin{cases}
(u_L,  \sigma_L),  & \text{if }  x <(u_L + k)t +x_0,  \forall t ,\\
(\frac{x - x_0}{t} - k, - k \frac{x - x_0}{t} + \sigma_L + k (u_L + k)), & \text{if }  (u_L + k)t +x_0 < x< (u_m + k)t  +x_0, \forall t,\\
(u_m,  \sigma_m ),   & \text{if  }  (u_m + k)t +x_0  <x< (u_m   + k)t  +x_1, \forall t, \\
(\frac{x - x_1}{t} - k, - k \frac{x - x_1}{t} + \sigma_m + k (u_m + k) ) ,    & \text{if  } (u_m + k)t +x_1<x< (u_R + k)t +x_1, \forall t, \\
(u_R,  \sigma_R),   & \text{if  } x >  (u_R + k)t +x_1,  \forall t,
\end{cases}
\label{4.5}
\end{gathered}
 \end{equation}
see figure 3. \\

Case 2:  When $(u_m, \sigma_m) \in S_2( (u_L, \sigma_L)  )$ \\
$u_L$ is connected to $u_m$ by a 2 - shock as seen in Figure 4 below. \\

\begin{figure}[!hbtp]
\includegraphics[width=12cm,height=4cm]{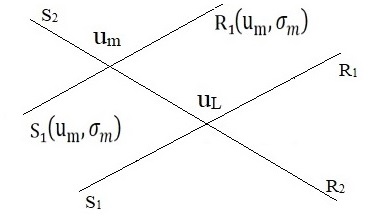}
\caption{ When $(u_m, \sigma_m) \in S_2( (u_L, \sigma_L)  )$}
\end{figure} 
Subcase 1:   If $(u_R, \sigma_R) \in R_2( (u_L, \sigma_L)  ) \bigcup \{  S_2 (u_L, \sigma_L), u_m< u_R< u_L   \}:$ \\
 then, $u_m$ is connected  to $u_R$ by a $2-$rarefaction. Now, since $u_L> u_m$ we have $\frac{u_L +u_m}{2} + k > u_m + k.$  

\begin{figure}[!hbtp]
\includegraphics[width=14cm,height=4cm]{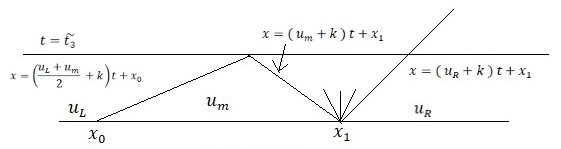}
\caption{ For Case 2,  Subcase 1}
\end{figure} 

So  the lines  $x= (\frac{u_L +u_m}{2} + k) t +x_0$ and $x= (u_m + k)t +x_1$ intersect at $( \tilde{x_3}, \tilde{t_3} ).$ The curve for $t > \tilde{t_3},$ is given by, 
\begin{equation}
\begin{aligned}
\frac{dx}{dt}& = \frac{ \frac{x -  x_1}{ t } + u_L - k}{2 }
 - \frac{ - k \frac{x - x_1}{t} + \sigma_m + k (u_m + k) - \sigma_L }{  \frac{x - x_1}{t} -k - u_L }, \\
\frac{dx}{dt} &= \frac{ \frac{x -  x_1}{ t } + u_L - k}{2 } - k^2 ( \frac{ \frac{x - x_1}{t} -k - u_L }{  - k \frac{x - x_1}{t} + \sigma_m + k (u_m + k) - \sigma_L }).
\label{4.6}
\end{aligned} 
 \end{equation}
Now $(u_m, \sigma_m ) \in S_2(u_L, \sigma_L )$  implies $\sigma_m - \sigma_L =  - k ( u_m - u_L ).$ This with \eqref{4.6} gives
\begin{equation*}
\frac{dx}{dt} = \frac{u_L + k }{2} +  \frac{x - x_1}{2 t}, \,\,\, x(\tilde{t_3}) = \tilde{x_3},
 \end{equation*}
which  in turn gives the curve  $x= \beta_2(t)$ for $t > \tilde{t_3},$ as
$x = (u_L + k) t + c t^{\frac{1}{2}} +x_0$ where $c = (\tilde{x_3} - x_0) {\tilde{t_3}}^{- \frac{1}{2}} - (u_L + k) {\tilde{t_3}}^{\frac{1}{2}}.$ 
Here there are two cases. The curves $x= \beta_2(t)$ and $x= ( u_R + k) t +x_1,$ may or may not intersect. If the two curves donot intersect then the solution is
\begin{equation}
\begin{gathered}
( u(x,t), \sigma(x, t) ) = \\
\begin{cases}
(u_L, \sigma_L ),  & \text{if } 0< x < (\frac{u_L +u_m}{2} + k) t +x_0, t \leq \tilde{t_3} ,\\
(u_m, \sigma_m ), & \text{if }  (\frac{u_L +u_m}{2} + k) t +x_0 <  x<(u_m + k)t  +x_1, 
\\ & \,\,\, \,\,\, \text{and } t \leq \tilde{t_3},\\
( \frac{x - x_1}{t} -k, - k \frac{x - x_1}{t} + \sigma_m + k (u_m + k) ),   & \text{if  } (u_m + k)t +x_1 <x < (u_R + k)t +x_1,  \\ & \,\,\, \,\,\, \text{and } t \leq \tilde{t_3}, \\    
  (u_R, \sigma_R),     & \text{if  }   x> (u_R + k)t +x_1, \forall t, \\
(u_L, \sigma_L ),  & \text{if } 0< x <\beta_2(t), t >  \tilde{t_3},\\
( \frac{x - x_1}{t} -k, - k \frac{x - x_1}{t} + \sigma_m + k (u_m + k) ),  & \text{if } \beta_2(t)< x< (u_R + k)t +x_1, t >  \tilde{t_3}.\\
\end{cases}
\label{4.7}
\end{gathered}
 \end{equation}
 If the curves $x= \beta_2(t)$ and $x= (u_R  + k) t +x_1,$  intersect at $(x_3,t_3)$ and if $(u_R, \sigma_R) \in R_2(u_L, \sigma_L),$  then for $t>t_3$, solution is a $2-$ rarefaction with left value $(u_L,\sigma_L)$ and right value $(u_R,\sigma_R)$. So   the solution then is,
\begin{equation}
\begin{gathered}
( u(x,t), \sigma(x, t) )= \\
 \begin{cases}
(u_L, \sigma_L ),  & \text{if } 0< x < (\frac{u_L +u_m}{2}  + k) t +x_0, t \leq \tilde{t_3} ,\\
(u_m, \sigma_m ), & \text{if }  (\frac{u_L +u_m}{2} + k) t +x_0 <   x<(u_m + k)t +x_1, 
\\ & \,\,\, \,\,\, \text{and } t \leq \tilde{t_3},\\
( \frac{x - x_1}{t} -k, - k \frac{x - x_1}{t} + \sigma_m + k (u_m + k) ),   & \text{if  } (u_m + k)t +x_1 <x  < (u_R + k)t +x_1,  \\ & \,\,\, \,\,\, \text{and } t \leq \tilde{t_3}, \\    
  (u_R, \sigma_R),     & \text{if  }   x> (u_R + k)t +x_1, \forall t,\\
(u_L, \sigma_L ), & \text{if }  0 < x<  \beta_2(t), \tilde{t_3}< t < t_3, \\
 ( \frac{x - x_1}{t} -k, - k \frac{x - x_1}{t} + \sigma_m + k (u_m + k) ),  & \text{if }   \beta_2(t)< x<  (u_R + k)t  +x_1,   \tilde{t_3}< t < t_3, \\
(u_L, \sigma_L ), & \text{if  } 0< x<  (u_L + k)(t - t_3) +x_3,  t > t_3, \\
( \frac{x - x_1}{t} -k, - k \frac{x - x_1}{t} + \sigma_L + k (u_L + k) ),   
& \text{if  }  (u_L + k))(t - t_3) +x_3 < x< (u_R + k))(t - t_3) 
\\ & \,\,\, \,\,\, +x_3,  t > t_3, \\
 (u_R, \sigma_R),     & \text{if  }  x>  (u_R + k))(t - t_3) +x_3,  t > t_3. \\
\end{cases}
\label{4.8}
\end{gathered}
 \end{equation}
 If the curves $x= \beta_2(t)$ and $x= (u_R  + k) t +x_1,$  intersect at $(x_3,t_3)$ and if $(u_R, \sigma_R) \in S_2(u_L, \sigma_L), u_m < u_R< u_L,$  then for $t>t_3$, solution is a $2-$ shock with left value $(u_L,\sigma_L)$ and right value $(u_R,\sigma_R)$. So the solution  then is,
\begin{equation}
\begin{gathered}
( u(x,t), \sigma(x, t) )= \\
 \begin{cases}
(u_L, \sigma_L ),  & \text{if } 0< x < (\frac{u_L +u_m}{2}  + k) t +x_0, t \leq \tilde{t_3} ,\\
(u_m, \sigma_m ), & \text{if }  (\frac{u_L +u_m}{2} + k) t +x_0 <   x<(u_m + k)t  +x_1, 
\\ & \,\,\, \,\,\, \text{and } t \leq \tilde{t_3},\\
( \frac{x - x_1}{t} -k, - k \frac{x - x_1}{t} + \sigma_m + k (u_m + k) ),   & \text{if  } (u_m + k)t +x_1 <x  < (u_R + k)t +x_1,  \\ & \,\,\, \,\,\, \text{and } t \leq \tilde{t_3}, \\    
  (u_R, \sigma_R),     & \text{if  }   x> (u_R + k)t +x_1, \forall t,\\
(u_L, \sigma_L ), & \text{if }  0 < x<  \beta_2(t), \tilde{t_3}< t < t_3, \\
 ( \frac{x - x_1}{t} -k, - k \frac{x - x_1}{t} + \sigma_m + k (u_m + k) ),  & \text{if }   \beta_2(t)< x<  (u_R + k)t  +x_1,  \tilde{t_3}< t < t_3, \\
(u_L, \sigma_L ), & \text{if  } 0< x<  (\frac{u_L + u_R}{2} + k))(t - t_3) +x_3,  t > t_3, \\
 (u_R, \sigma_R),     & \text{if  }  x>  (\frac{u_L + u_R}{2} + k))(t - t_3) +x_3,  t > t_3, \\
\end{cases}
\label{4.9}
\end{gathered}
 \end{equation}
see figure 5. \\

Subcase 2:   If $(u_R, \sigma_R) \in S_2( (u_L, \sigma_L)  ), \,\,\,  u_R < u_m < u_L:$ \\ 
 then, $u_m$ is connected  to $u_R$ by a $2-$rarefaction. Now, since $u_L> u_m$ we have $\frac{u_L +u_m}{2} + k > u_m + k.$  

\begin{figure}[!hbtp]
\includegraphics[width=14cm,height=4cm]{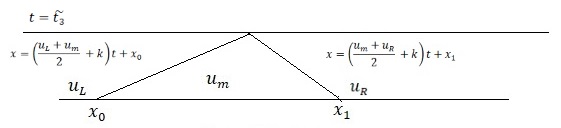}
\caption{ For Case 2,  Subcase 2}
\end{figure} 
So  the lines  $x= (\frac{u_L +u_m}{2} + k) t +x_0$ and $x= (\frac{u_R +u_m}{2} + k) t  +x_1$ intersect at $( \tilde{x_3}, \tilde{t_3} ).$ The curve for $t > \tilde{t_3}$ is given by, 
\begin{equation}
\begin{aligned}
\frac{dx}{dt}& = \frac{ \frac{x -  x_1}{ t } + u_L - k}{2 }
 - \frac{ - k \frac{x - x_1}{t} + \sigma_m + k (u_m + k) - \sigma_L }{  \frac{x - x_1}{t} -k - u_m }, \\
\frac{dx}{dt} &= \frac{ \frac{x -  x_1}{ t } + u_L - k}{2 } - k^2 ( \frac{ \frac{x - x_1}{t} -k - u_m }{  - k \frac{x - x_1}{t} + \sigma_m + k (u_m + k) - \sigma_L }).
\label{4.10}
\end{aligned} 
 \end{equation}
Using $(u_m, \sigma_m ) \in S_2(u_L, \sigma_L )$ that implies $\sigma_m - \sigma_L =  - k ( u_m - u_L )$ in \eqref{4.10} we get,
\begin{equation*}
\frac{dx}{dt} = \frac{u_L + k }{2} +  \frac{x - x_1}{2 t}, \,\,\, x(\tilde{t_3}) = \tilde{x_3},
 \end{equation*}
which  in turn gives the curve  $x= \beta_3(t)$ for $t > \tilde{t_3},$ as
$x = (u_L + k) t + c t^{\frac{1}{2}} +x_0$ where $c = (\tilde{x_3} - x_0) {\tilde{t_3}}^{- \frac{1}{2}} - (u_L + k) {\tilde{t_3}}^{\frac{1}{2}}.$ 
 The solution then is
\begin{equation}
\begin{gathered}
( u(x,t), \sigma(x, t) )= \begin{cases}
(u_L, \sigma_L ),  & \text{if } 0< x < (\frac{u_L +u_m}{2} + k) t   +x_0, t \leq \tilde{t_3} ,\\
(u_m, \sigma_m), & \text{if }  (\frac{u_L +u_m}{2} + k) t +x_0 < x  < (\frac{u_R +u_m}{2} + k) t  +x_1, t \leq \tilde{t_3},\\
(u_R, \sigma_R), & \text{if } x> (\frac{u_R +u_m}{2} + k) t  +x_1, \forall t, \\
 (u_L, \sigma_L ),  & \text{if } 0 < x<   \beta_3(t),    t > \tilde{t_3}, \\
 (u_R, \sigma_R),  & \text{if } x> \beta_3(t),  t > \tilde{t_3},
\end{cases}
\label{4.11}
\end{gathered}
 \end{equation}
see figure 6.\\ 

 Next we consider $u_m, u_L, u_R$ to lie on the level set of the first Riemann Invariant. \\

Case 3: When $(u_m, \sigma_m) \in S_1( (u_L, \sigma_L)  )$ \\
$u_L$ is connected to $u_m$ by a $1$-shock as seen in Figure 7 below. \\
\begin{figure}[!hbtp]
\includegraphics[width=12cm,height=3cm]{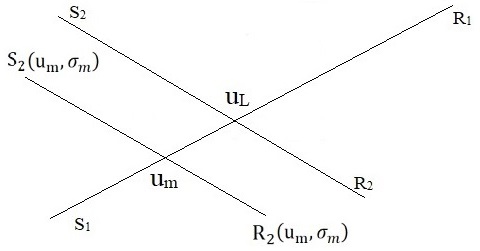}
\caption{ When $(u_m, \sigma_m) \in S_1( (u_L, \sigma_L)  )$}
\end{figure} 

Subcase 1: If $( u_R, \sigma_R ) \in R_1(u_L, \sigma_L) \bigcup \{ S_1(u_L, \sigma_L), u_m <u_R<u_L \}:$\\
$u_m$ is connected to $u_R$ by a $1$-rarefaction. 

\begin{figure}[!hbtp]
\includegraphics[width=14cm,height=4cm]{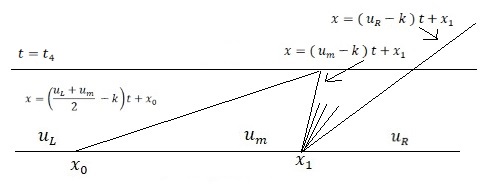}
\caption{ For Case 3,  Subcase 1}
\end{figure} 
Now since $u_L> u_m$ we have $\frac{u_L +u_m}{2} - k > u_m - k.$  So  the lines  $x= (\frac{u_L +u_m}{2} - k) t +x_0$ and $x= (u_m - k) t +x_1$ intersect at $( x_4, t_4 ).$ The curve for $t >t_4$ is given by, 
\begin{equation}
\begin{aligned}
\frac{dx}{dt} &= \frac{u_L + \frac{x - x_1}{t} +k}{2} - \frac{  k \frac{x - x_1}{t} + \sigma_m - k (u_m - k)  - \sigma_L}{ \frac{x - x_1}{t} +k - u_L }, \\
\frac{dx}{dt} &= \frac{u_L + \frac{x - x_1}{t} +k}{2} - k^2 ( \frac{  \frac{x - x_1}{t} +k - u_L  }{  k \frac{x - x_1}{t} + \sigma_m - k (u_m - k)  - \sigma_L }).
\end{aligned}
\label{4.12}
 \end{equation}
Using $(u_m, \sigma_m) \in S_1(u_L, \sigma_L) $  that gives
$\sigma_L - \sigma_m = k (u_L - u_m )$ with  \eqref{4.12} we get
\begin{equation*}
\frac{dx}{dt} = \frac{u_L - k }{2} +  \frac{x - x_1}{2 t}, \,\,\, x(t_4) = x_4,
 \end{equation*}
which in turn gives the curve  $x= \phi_1(t)$ for $t > t_4,$ as
$ \phi_1(t)= (u_L - k) t + c t^{\frac{1}{2}} +x_1$ where $c=  (x_4 - x_0) {t_4}^{- \frac{1}{2}} - (u_L - k) {t_4}^{\frac{1}{2}}.$ The curves $x= \phi_1(t)$ and $x= (u_R - k) t +x_1,$  may or may not intersect. If these curves donot intersect, the solution is,
\begin{equation}
\begin{gathered}
( u(x, t), \sigma(x, t) ) \\
= \begin{cases}
(u_L, \sigma_L ),  & \text{if } 0< x < (\frac{u_L +u_m}{2} - k) t   +x_0,   t \leq t_4,\\
(u_m, \sigma_m), & \text{if }  (\frac{u_L +u_m}{2} - k) t +x_0 < x  <  (u_m  - k) t +x_1, 
\\ & \,\,\, \,\,\, \text{and } t \leq t_4,\\
(\frac{x - x_1}{t} +k,   k \frac{x - x_1}{t} + \sigma_m - k (u_m - k) ), & \text{if }  (u_m - k) t +x_1 < x<   (u_R - k) t +x_1, t \leq t_4,\\
(u_R, \sigma_R), & \text{if }   x>  (u_R - k) t +x_1, \forall t, \\
(u_L, \sigma_L ) & \text{if }  0 < x <  \phi_1(t),  t > t_4,\\
(\frac{x - x_1}{t} +k,   k \frac{x - x_1}{t} + \sigma_m - k (u_m - k) ), & \text{if } \phi_1(t)<x<  (u_R - k) t  +x_1, t > t_4.
\end{cases}
\label{4.13}
\end{gathered}
 \end{equation}
 If the curves $x= \phi_1(t)$ and $x= (u_R  + k) t +x_1$  intersect at $(x_3,t_3)$ and if $(u_R, \sigma_R) \in R_1(u_L, \sigma_L),$  then for $t>t_3$,  we have a Riemann problem's solution which  is a $1-$ rarefaction with left value $(u_L,\sigma_L)$ and right value $(u_R,\sigma_R)$. So the solution then is,
\begin{equation}
\begin{gathered}
( u(x, t), \sigma(x, t) ) \\
= \begin{cases}
(u_L, \sigma_L ),  & \text{if } 0< x < (\frac{u_L +u_m}{2} - k) t  +x_0,  t \leq t_4,\\
(u_m, \sigma_m), & \text{if }  (\frac{u_L +u_m}{2} - k) t +x_0 < x  <  (u_m - k) t +x_1, 
\\ & \,\,\, \,\,\, \text{and } t \leq t_4,\\
(\frac{x - x_1}{t} +k,   k \frac{x - x_1}{t} + \sigma_m - k (u_m - k) ), & \text{if }  (u_m - k) t +x_1 < x<  (u_R - k) t +x_1,
\\ & \,\,\, \,\,\, \text{and } t \leq t_4,\\
(u_R, \sigma_R), & \text{if }   x>  (u_R - k) t +x_1, \forall t, \\
(u_L, \sigma_L ) & \text{if }  0 < x <  \phi_1(t),   t_4< t< t_3,\\
(\frac{x - x_1}{t} +k,   k \frac{x - x_1}{t} + \sigma_m - k (u_m - k) ), & \text{if } \phi_1(t)<x<  (u_R - k) t   +x_1,  t_4< t< t_3, \\
(u_L, \sigma_L ),  & \text{if }  0 < x<  (u_L - k) (t - t_3) +x_3,   t> t_3, \\
(\frac{x - x_1}{t} +k,   k \frac{x - x_1}{t} + \sigma_m - k (u_m - k) ), & \text{if }  (u_L - k)  (t - t_3) +x_3< x< (u_R - k)  (t - t_3) 
\\ & \,\,\, \,\,\,  +x_3,    t> t_3, \\
(u_R, \sigma_R), & \text{if }  x>  (u_R - k)  (t - t_3)+x_3, t> t_3.
\end{cases}
\label{4.14}
\end{gathered}
 \end{equation}
 If the curves $x= \phi_1(t)$ and $x= (u_R  + k) t +x_1,$  intersect at $(x_3,t_3)$ and if $(u_R, \sigma_R) \in S_1(u_L, \sigma_L),$ with $u_m <u_R<u_L$ then for $t>t_3$, we have a Riemann problem's solution which is a $1-$ shock with left value $(u_L,\sigma_L)$ and right value $(u_R,\sigma_R)$. So  the solution then is,
\begin{equation}
\begin{gathered}
( u(x, t), \sigma(x, t) ) \\
= \begin{cases}
(u_L, \sigma_L ),  & \text{if } 0< x < (\frac{u_L +u_m}{2} - k) t  +x_0,   t \leq t_4,\\
(u_m, \sigma_m), & \text{if }  (\frac{u_L +u_m}{2} - k) t +x_0 < x  <  (u_m  - k) t +x_1,
\\ & \,\,\, \,\,\, \text{and } t \leq t_4,\\
(\frac{x - x_1}{t} +k,   k \frac{x - x_1}{t} + \sigma_m - k (u_m - k) ), & \text{if }  (u_m - k) t +x_1 < x<  (u_R - k) t  +x_1, t \leq t_4,\\
(u_R, \sigma_R), & \text{if }   x>  (u_R - k) t +x_1, \forall t, \\
(u_L, \sigma_L ) & \text{if }  0 < x <  \phi_1(t),   t_4< t< t_3,\\
(\frac{x - x_1}{t} +k,   k \frac{x - x_1}{t} + \sigma_m - k (u_m - k) ), & \text{if } \phi_1(t)<x<  (u_R - k) t   +x_1,  t_4< t< t_3, \\
(u_L, \sigma_L ),  & \text{if }  0 < x<   (\frac{u_L +u_R}{2} - k)  (t - t_3) +x_3,   t> t_3, \\
(u_R, \sigma_R), & \text{if }  x> (\frac{u_L +u_R}{2} - k)  (t - t_3) +x_3, t> t_3,
\end{cases}
\label{4.15}
\end{gathered}
 \end{equation}
see figure 8. \\

Subcase 2: If $( u_R, \sigma_R ) \in S_1((u_L, \sigma_L)), \,\,\, u_R < u_m < u_L,$ \\ 
$u_m$ is connected to $u_R$ by a $1$-shock. Now since $u_L> u_R$ we have $\frac{u_L +u_m}{2} - k > \frac{u_R +u_m}{2} - k.$ 
\begin{figure}[!hbtp]
\includegraphics[width=14cm,height=4cm]{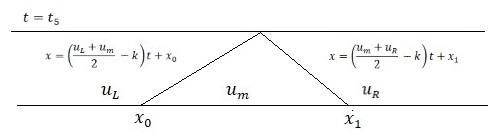}
\caption{ For Case 3,  Subcase 2}
\end{figure} 

 So  the lines  $x= (\frac{u_L +u_m}{2} - k) t +x_0$ and $x= ( \frac{u_R +u_m}{2} - k) t +x_1$ intersect at $( x_5, t_5 ).$ The curve for $t > t_5$ is given by, 
\begin{equation}
\begin{aligned}
\frac{dx}{dt} &=( \frac{u_L + u_R}{2 } ) - \frac{ \sigma_R - \sigma_L}{u_R - u_L}, \\
  \frac{dx}{dt} &=( \frac{u_L + u_R}{2 } ) - k^2 ( \frac{u_R - u_L}{\sigma_R - \sigma_L}), 
\label{4.16}
\end{aligned}
 \end{equation}
$( u_R, \sigma_R ) \in S_1(u_L, \sigma_L)$ gives $ \sigma_R - \sigma_L = k ( u_R - u_L )$ which with \eqref{4.16} gives  the $\frac{dx}{dt} =( \frac{u_L + u_R}{2 } )  - k $  with $x(t_5) = x_5. $ This gives the curve  $x = ( \frac{u_L + u_R}{2 }   - k ) ( t -  t_5)+ x_5$  for  $t > t_5.$  The solution for this case is,
\begin{equation}
\begin{gathered}
( u(x, t), \sigma(x, t) )\\
= \begin{cases}
(u_L, \sigma_L),  & \text{if } 0< x < (\frac{u_L +u_m}{2} - k) t +x_0, t \leq t_5,\\
(u_m, \sigma_m), & \text{if }  (\frac{u_L +u_m}{2} - k) t +x_0 < x< (\frac{u_R +u_m}{2} - k) t +x_1, t \leq t_5,\\
(u_R, \sigma_R), & \text{if }  x>  (\frac{u_R +u_m}{2} - k) t +x_1, t \leq t_5,\\
(u_L, \sigma_L), & \text{if }  0 < x<  ( \frac{u_L + u_R}{2 }   - k ) ( t -  t_5 )+ x_5, t > t_5, \\
 (u_R, \sigma_R), & \text{if } x > ( \frac{u_L + u_R}{2 }   - k ) ( t -  t_5 )+ x_5, t >t_5.
\end{cases}
\label{4.17}
\end{gathered}
\end{equation}
see figure 9. \\

Case 4: When $(u_m, \sigma_m) \in R_1( (u_L, \sigma_L)  )$ \\
$u_L$ is connected to $u_m$ by a $1$-rarefaction as seen in Figure 10 below. \\
\begin{figure}[!hbtp]
\includegraphics[width=12cm,height=4cm]{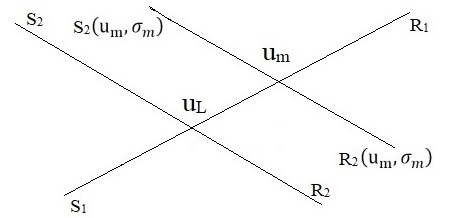}
\caption{ When $(u_m, \sigma_m) \in R_1( (u_L, \sigma_L)  )$}
\end{figure}

Subcase 1: If $( u_R, \sigma_R ) \in S_1((u_L, \sigma_L)) \bigcup \{  R_1((u_L, \sigma_L)), u_L< u_R <u_m  \}$  \\
$u_m$ is connected to $u_R$ by a $1$-shock. Now since $u_m> u_R,$ we have $u_m - k > \frac{u_R +u_m}{2} - k.$  So  the lines  $x= (\frac{u_R +u_m}{2} - k) t +x_1$ and $x= (u_m - k) t +x_0$ intersect at $( x_6, t_6 ).$ The curve for $t > t_6$ is given by, 
\begin{figure}[!hbtp]
\includegraphics[width=14cm,height=4cm]{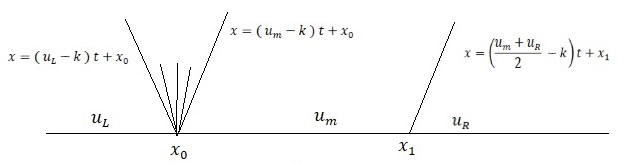}
\caption{ For Case 4,  Subcase 1}
\end{figure} 
\begin{equation}
\begin{aligned}
\frac{dx}{dt} &= \frac{ u_R +  \frac{x - x_0}{t} + k  }{2} - \frac{\sigma_R - (  k \frac{x - x_0}{t} + \sigma_L - k (u_L - k )}{ u_R -  ( \frac{x - x_0}{t} + k ) }, \\
\frac{dx}{dt} &= \frac{ u_R +  \frac{x - x_0}{t} + k  }{2} - k^2 ( \frac{ u_R -  ( \frac{x - x_0}{t} + k ) }{\sigma_R - (  k \frac{x - x_0}{t} + \sigma_L - k (u_L - k )} ),
\label{4.18}
\end{aligned}
 \end{equation}
Using $ \sigma_R - \sigma_L = k ( u_R - u_L )$ that comes from $(u_R, \sigma_R ) \in S_1(u_L, \sigma_L ) $  in \eqref{4.18}, we get  $\frac{dx}{dt} = \frac{ u_R +  \frac{x - x_0}{t} + k  }{2} - k $ which gives the curve $x = \phi_2(t)$ for $ t > t_6$ with  $\phi_2(t) =  (u_R - k) t + c t^{\frac{1}{2}} +x_1$ where $c=  (x_6 - x_0) {t_6}^{- \frac{1}{2}} - (u_R - k) {t_6}^{\frac{1}{2}}.$  The curves $x= \phi_2(t)$ and $x= (u_L - k) t +x_0$  may or may not intersect. If these curves donot intersect, the solution is,
\begin{equation}
\begin{gathered}
( u(x, t), \sigma(x, t) ) \\
= \begin{cases}
( u_L, \sigma_L),  & \text{if } 0< x < ( u_L - k) t +x_0, \forall t,\\
( \frac{x - x_0}{t} + k,  k \frac{x - x_0}{t} + \sigma_L - k (u_L - k) ), & \text{if }  ( u_L - k) t +x_0 < x< ( u_m  - k) t  +x_0,   t \leq t_6,\\
( u_m, \sigma_m), & \text{if }   ( u_m - k) t +x_0 < x<   (\frac{u_m +u_R}{2}  - k) t  +x_1,  
\\ & \,\,\, \,\,\, \text{and }  t \leq t_6,\\
( u_R, \sigma_R), & \text{if }    x>    (\frac{u_m +u_R}{2} - k) t  +x_1,  \forall t,\\
( \frac{x - x_0}{t} + k,  k \frac{x - x_0}{t} + \sigma_L - k (u_L - k) ), & \text{if } ( u_L - k) t +x_0< x< \phi_2(t),  t > t_6, \\
( u_R, \sigma_R), & \text{if } x> \phi_2(t), t > t_6.
\end{cases}
\label{4.19}
\end{gathered}
 \end{equation}
 If the curves $x= \phi_2(t)$ and $x= (u_L - k) t +x_0,$  intersect at $(x_7,t_7)$ and if $(u_R, \sigma_R) \in S_1(u_L, \sigma_L),$  then for $t>t_7$,  we have a Riemann problem's solution which  is a $1-$ shock with left value $(u_L,\sigma_L)$ and right value $(u_R,\sigma_R)$. So  the solution then is, \\
\begin{equation}
\begin{gathered}
( u(x, t), \sigma(x, t) ) \\
= \begin{cases}
( u_L, \sigma_L),  & \text{if } 0< x < ( u_L - k) t +x_0, \forall t,\\
( \frac{x - x_0}{t} + k,  k \frac{x - x_0}{t} + \sigma_L - k (u_L - k) ), & \text{if }  ( u_L - k) t +x_0 < x< ( u_m - k) t +x_0, 
\\ & \,\,\, \,\,\, \text{and }  t \leq t_6,\\
( u_m, \sigma_m), & \text{if }   ( u_m - k) t +x_0 < x<   (\frac{u_m +u_R}{2}  - k) t  +x_1, 
\\ & \,\,\, \,\,\, \text{and }  t \leq t_6,\\
( u_R, \sigma_R), & \text{if }    x>    (\frac{u_m +u_R}{2} - k) t  +x_1,  \forall t,\\
( \frac{x - x_0}{t} + k,  k \frac{x - x_0}{t} + \sigma_L - k (u_L - k) ), & \text{if } ( u_L - k) t +x_0< x< \phi_2(t),  t_6< t< t_7, \\
( u_R, \sigma_R), & \text{if } x> \phi_2(t),  t_6< t< t_7,\\
( u_L, \sigma_L),  & \text{if } 0< x <  (\frac{u_L +u_R}{2} - k) (t - t_7)  +x_7, t>t_7, \\
( u_R, \sigma_R),  & \text{if } x >  (\frac{u_L +u_R}{2} - k) (t - t_7)  +x_7,  t>t_7. \\
\end{cases}
\label{4.20}
\end{gathered}
 \end{equation}

 If the curves $x= \phi_2(t)$ and $x= (u_L  - k) t +x_0,$  intersect at $(x_7,t_7)$ and if $(u_R, \sigma_R) \in R_1(u_L, \sigma_L),$ with $u_m <u_R<u_L$ then for $t>t_7$, we have a Riemann problem's solution which is a $1-$ rarefaction with left value $(u_L,\sigma_L)$ and right value $(u_R,\sigma_R)$. So  the solution then is,

\begin{equation}
\begin{gathered}
( u(x, t), \sigma(x, t) ) \\
= \begin{cases}
( u_L, \sigma_L),  & \text{if } 0< x < ( u_L - k) t +x_0, \forall t,\\
( \frac{x - x_0}{t} + k,  k \frac{x - x_0}{t} + \sigma_L - k (u_L - k) ), & \text{if }  ( u_L - k) t +x_0 < x< ( u_m  - k) t   +x_0,   t \leq t_6,\\
( u_m, \sigma_m), & \text{if }   ( u_m - k) t +x_0 < x<   (\frac{u_m +u_R}{2}   - k) t  +x_1, 
\\ & \,\,\, \,\,\, \text{and }  t \leq t_6,\\
( u_R, \sigma_R), & \text{if }    x>    (\frac{u_m +u_R}{2} - k) t  +x_1,  \forall t,\\
( \frac{x - x_0}{t} + k,  k \frac{x - x_0}{t} + \sigma_L - k (u_L - k) ), & \text{if } ( u_L - k) t +x_0< x< \phi_2(t),    t_6< t< t_7, \\
( u_R, \sigma_R), & \text{if } x> \phi_2(t),  t_6< t< t_7,\\
( u_L, \sigma_L),  & \text{if } 0< x <  (u_L - k) (t - t_7)  +x_7,   t>t_7, \\
( \frac{x - x_7}{t} + k,  k \frac{x - x_7}{t} + \sigma_L - k (u_L - k) ),  & \text{if }   (u_L - k) (t - t_7)  +x_7< x<  (u_R - k) (t - t_7) 
\\ & \,\,\, \,\,\,  +x_7 , t>t_7, \\
( u_R, \sigma_R),  & \text{if } x >  (u_R - k) (t - t_7)  +x_7,  t>t_7. \\
\end{cases}
\label{4.21}
\end{gathered}
 \end{equation}
see figure 11. \\

Subcase 2: If $( u_R, \sigma_R ) \in R_1((u_L, \sigma_L)),  \,\,\, u_L<u_m<u_R,$ \\ 
\begin{figure}[!hbtp]
\includegraphics[width=14cm,height=4cm]{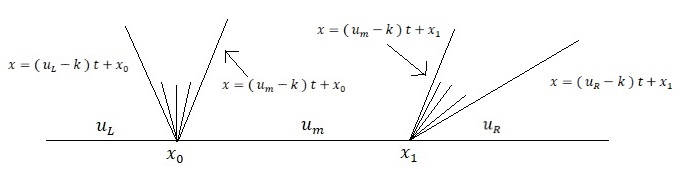}
\caption{ For Case 4,  Subcase 2}
\end{figure} 
$u_m$ is connected to $u_R$ by a $1$-rarefaction.
\begin{equation}
\begin{gathered}
( u(x, t), \sigma(x, t) )= \\
 \begin{cases}
(u_L, \sigma_L ),  & \text{if } 0< x < ( u_L - k) t +x_0, \forall t,\\
(\frac{x - x_0}{t} + k,  k \frac{x - x_0}{t} + \sigma_L - k (u_L - k) ), & \text{if }  ( u_L - k) t +x_0 < x< ( u_m - k) t  +x_0, \forall t,\\
(u_m, \sigma_m), & \text{if }   ( u_m - k) t +x_0 < x<   ( u_m - k) t   +x_1, \forall t,\\
( \frac{x - x_1}{t} + k,  k \frac{x - x_1}{t} + \sigma_m - k (u_m - k) ), & \text{if }   ( u_m - k) t +x_0 < x<   ( u_R - k) t  +x_1, \forall t,\\
(u_R, \sigma_R), & \text{if }    x>    ( u_R - k) t +x_1,  \forall t,
\end{cases}
\label{4.22}
\end{gathered}
 \end{equation}
see figure 12.  This elaborates the cases when  $( u_L, \sigma_L ), ( u_m, \sigma_m ), ( u_R, \sigma_R ),$ all lie on the level set of the same Riemann Invariants.  Explicit formulae for the solution given in \eqref{4.2} - \eqref{4.5}, \eqref{4.7} - \eqref{4.9}, \eqref{4.11}, \eqref{4.13} - \eqref{4.15}, \eqref{4.17} and \eqref{4.19} - \eqref{4.22} prove the theorem. \\

By construction we have shown that the formula given for $(u, \sigma)$ satisfies the equation in the region of smoothness and  the Rankine - Hugoniot condition \eqref{m2.7} along the curve of discontinuity holds. So $(u, \sigma)$ satisfies \eqref{e1.1} in the sense of measures \eqref{m2.4} with initial condition \eqref{e1.2}. \\

{ \bf Remark:} The cases when they do not all lie on the same  Riemann Invariants is more complicated due to there  being, after a certain finite time t,  either  interaction between 2 rarefactions in some cases or in some other cases interaction between a rarefaction and a shock, which is difficult to handle. These cases together with the cases explored in the above  theorem  cover all possible cases, when the three states of the initial data lie in any part of the domain. The analysis for these regions being more complicated is part of ongoing work as a continuation of this work. The novelty of the work is in not assuming any smallness  on  $( u_L, \sigma_L ), ( u_m, \sigma_m ), ( u_R, \sigma_R )$. 

\section{Conclusion}
 This paper is concerned with the interaction of Riemann solution centered about two points. Understanding the Riemann problem and  study of interaction of waves originating from two different points are important
steps in the construction of solution using approximation procedure like the Glimm's scheme. This amounts to constructing solution  with initial data consisting of three states $(u_L, \sigma_L), (u_m, \sigma_m),$ and  $(u_R, \sigma_R)$. We solve the system without any smallness condition of the data for certain class of states. Since the system is nonconservative we need to choose the nonconservative product in addition to admisibility condition for the shock. We use the Volpert product \cite{v1} and Lax entropy inequality for the shocks. \\

The problem considered in this work, for arbitrary data is still an open problem and has a quite complicated structure.
For general  non-conservative systems, most previous existing results about  the Riemann problem and their interactions  require smallness conditions on the data for the analysis. We use the special structure of the Elastodynamics system to provide an analysis that  does not need this smallness conditions on the data, with the assumption that all data lie in the level set of one of the Riemann invariants.

\section*{Acknowledgments}

Funding: No funding involved. \\
Conflicts of interest/Competing interests: None. \\
Availability of data and material: Not applicable. \\
Code availability: Not applicable.

\section*{ORCID}
Kayyunnapara Divya Joseph https://orcid.org/0000-0002-4126-7882

\end{document}